\documentclass[12pt,reqno]{amsart}
\usepackage{enumerate, latexsym, amsmath, amsfonts, amssymb, amsthm, color}
\def\Z{\Bbb Z}
\def\N{\Bbb N}

\def\l{\left}
\def\r{\right}
\def\bg{\bigg}
\def\({\bg(}
\def\){\bg)}
\def\t{\text}
\def\f{\frac}

\def\ls{\leqslant}
\def\gs{\geqslant}

\def\bi{\binom}
\def\al{\alpha}

\def\Proof{\noindent{\it Proof}}

\theoremstyle{plain}
\newtheorem{theorem}{Theorem}

\newtheorem{lemma}{Lemma}

\theoremstyle{definition}

\theoremstyle{remark}
\newtheorem{remark}{Remark}

 \vspace{4mm}

\begin{document}
 \baselineskip=17pt
\hbox{Publ. Math. Debrecen 85(2014), no.\,3-4, 285--295.}
\medskip

\title
[On monotonicity of some combinatorial sequences] {On monotonicity
of some combinatorial sequences}

\author
[Qing-Hu Hou, Zhi-Wei Sun and Hao-Min Wen] {Qing-Hu Hou*, Zhi-Wei
Sun** and Haomin Wen}

\thanks{*Supported by the National Natural Science Foundation (grant 11171167) of China}

\thanks{**Supported by the National Natural Science Foundation (grant 11171140)
 of China}

\address {(Qing-Hu Hou) Center for Applied Mathematics, Tianjin University, Tianjin 300072, People's Republic of China}
\email{hou@nankai.edu.cn}

\address {(Zhi-Wei Sun) Department of Mathematics, Nanjing
University, Nanjing 210093, People's Republic of China}
\email{zwsun@nju.edu.cn}

\address{(Haomin Wen) Department of Mathematics
\\The University of Pennsyvania
\\Philadelphia, PA  19104, USA} \email{weh@math.upenn.edu}

\keywords{Combinatorial sequences, monotonicity, log-concavity
\newline \indent 2010 {\it Mathematics Subject Classification}. Primary 05A10; Secondary 11B39, 11B75.}

 \begin{abstract} We confirm Sun's conjecture that $(\root{n+1}\of{F_{n+1}}/\root{n}\of{F_n})_{n\gs 4}$
 is strictly decreasing to the limit 1, where $(F_n)_{n\gs0}$ is the Fibonacci sequence.
 We also prove that the sequence $(\root{n+1}\of{D_{n+1}}/\root{n}\of{D_n})_{n\gs3}$ is strictly decreasing with limit $1$,
 where $D_n$ is the $n$-th derangement number. For $m$-th order harmonic numbers $H_n^{(m)}=\sum_{k=1}^n 1/k^m\ (n=1,2,3,\ldots)$, we show that
$(\root{n+1}\of{H^{(m)}_{n+1}}/\root{n}\of{H^{(m)}_n})_{n\gs3}$ is
strictly increasing.
\end{abstract}

\maketitle

\section{Introduction}
\setcounter{lemma}{0}
\setcounter{theorem}{0}
\setcounter{corollary}{0}
\setcounter{remark}{0}
\setcounter{conjecture}{0}

A challenging conjecture of Firoozbakht states that
$$\root n\of{p_n}>\root{n+1}\of {p_{n+1}}\quad\mbox{for every}\ n=1,2,3,\ldots,$$
where $p_n$ denotes the $n$-th prime. Note that $\lim_{n\to\infty}\root n\of{p_n}=1$ by the Prime Number Theorem.
In \cite{S13} the second author conjectured further that
for any integer $n>4$ we have the inequality
$$\f{\root{n+1}\of{p_{n+1}}}{\root n\of{p_n}}<1-\f{\log\log n}{2n^2},$$
which has been verified for all $n\ls3.5\times10^6$.
Motivated by this and \cite{S13a}, Sun \cite[Conj. 2.12]{S13} conjectured that the sequence
$(\root{n+1}\of{S_{n+1}}/\root n\of{S_n})_{n\gs7}$ is strictly increasing,
where $S_n$ is the sum of the first $n$ positive squarefree numbers.
Moreover, he also posed many conjectures on
monotonicity of sequences of the type $(\root{n+1}\of{a_{n+1}}/\root
n\of{a_n})_{n\gs N}$ with $(a_n)_{n\gs1}$ a familiar combinatorial
sequence of positive integers.

Throughout this paper, we set $\N=\{0,1,2,\ldots\}$ and
$\Z^+=\{1,2,3,\ldots\}$.

Let $A$ and $B$ be integers with $\Delta=A^2-4B\not=0$. The Lucas sequence $u_n=u_n(A,B)\
(n\in\N)$ is defined as follows:
$$u_0=0,\ u_1=1,\ \mbox{and}\ u_{n+1}=Au_n-Bu_{n-1}\ \t{for}\
n=1,2,3,\ldots.$$
It is well known that $u_n=(\al^n-\beta^n)/(\al-\beta)$ for all $n\in\N$, where
$$\al=\f{A+\sqrt{\Delta}}2\  \ \mbox{and}\ \ \beta=\f{A-\sqrt{\Delta}}2$$
are the two roots of the characteristic equation $x^2-Ax+B=0$. The
sequence $F_n=u_n(1,-1)\ (n\in\N)$ is the famous Fibonacci sequence,
see \cite[p.\,46]{St} for combinatorial interpretations of Fibonacci
numbers.

Our first result is as follows.

\begin{theorem} \label{Th1.1} Let $A>0$ and $B\not=0$ be integers with $\Delta=A^2-4B>0$, and set $u_n=u_n(A,B)$ for $n\in\N$.
 Then there exists an integer $N>0$ such that the sequence
$(\root{n+1}\of{u_{n+1}}/\root n\of{u_n})_{n\gs N}$ is strictly
decreasing with limit $1$. In the case $A=1$ and $B=-1$ we may take
$N=4$.
\end{theorem}

\begin{remark}\label{Rem1.1} Under the condition of Theorem 1.1, by \cite[Lemma 4]{S92} we have $u_n<u_{n+1}$ unless $A=n=1$.
 Note that the second assertion in Theorem \ref{Th1.1} confirms a conjecture
 of the second author \cite[Conj. 3.1]{S13} on the Fibonacci sequence.
\end{remark}

For $n\in\Z^+$ the $n$th derangement number $D_n$ denotes the number
of permutations $\sigma$ of $\{1,\ldots,n\}$ with $\sigma(i)=i$ for
no $i=1,\ldots,k$. It has the following explicit expression (cf.
\cite[p.\,67]{St}):
$$D_n=\sum_{k=0}^n(-1)^k\f{n!}{k!}.$$

Our second theorem is the following result conjectured by the second author \cite[Conj. 3.3]{S13}.

\begin{theorem}\label{Th1.2} The sequence
    $\left(\sqrt[n+1]{D_{n+1}}/\sqrt[n]{D_n}\right)_{n \gs3}$ is strictly decreasing with limit $1$.
    \label{}
\end{theorem}
\begin{remark} It follows from Theorem 1.2 that the sequence $(\root n\of{D_n})_{n\gs2}$ is strictly increasing.
\end{remark}

For each $m\in\Z^+$ those $H_n^{(m)}=\sum_{k=1}^n1/k^m\ (n\in\Z^+)$ are called harmonic numbers of order $m$.
The usual harmonic numbers (of order 1)
are those rational numbers $H_n=H_n^{(1)}\ (n=1,2,3,\ldots)$.

Our following theorem confirms Conjecture 2.16 of Sun \cite{S13}.

\begin{theorem}\label{Th1.3} For any positive integer $m$, the sequence
    $\big(\sqrt[n+1]{H_{n+1}^{(m)}}/\sqrt[n]{H_n^{(m)}}\big)_{n \gs3}$ is strictly increasing.
    \label{}
\end{theorem}

We will prove Theorems 1.1--1.3 in Sections 2--4 respectively.
It seems that there is no simple form for the generating function $\sum_{n=0}^\infty\root n\of{a_n}x^n$
with $a_n=u_n, D_n, H_n^{(m)}$. Note also that the set of those sequences $(a_n)_{n\gs1}$ of positive numbers with $(\sqrt[n+1]{a_{n+1}}/\sqrt[n]{a_n})_{n\gs1}$
decreasing (or increasing) is closed under multiplication.

\section{Proof of Theorem 1.1}
\setcounter{lemma}{0}
\setcounter{theorem}{0}
\setcounter{corollary}{0}
\setcounter{remark}{0}
\setcounter{conjecture}{0}

\noindent{\it Proof of Theorem 1.1}. Set
$$\al=\f{A+\sqrt{\Delta}}2,\ \beta=\f{A-\sqrt{\Delta}}2,\ \mbox{and}\ \gamma=\f{\beta}{\al}=\f{A-\sqrt{\Delta}}{A+\sqrt{\Delta}}.$$
Then
$$\log u_n=\log\f{\al^n(1-\gamma^n)}{\al-\beta}=n\log\al+\log(1-\gamma^n)-\log\sqrt{\Delta}$$ for any $n\in\Z^+$.
Note that
$$\log\f{\root{n+1}\of{u_{n+1}}}{\root n\of{u_n}}=\f{\log u_{n+1}}{n+1}-\f{\log u_n}n
=\f{\log(1-\gamma^{n+1})}{n+1}-\f{\log(1-\gamma^n)}n+\f{\log\sqrt{\Delta}}{n(n+1)}.$$
Since
$$\lim_{n\to\infty}\f{\log(1-\gamma^n)}n=\lim_{n\to\infty}\f{-\gamma^n}n=0\ \mbox{and}\ \ \lim_{n\to\infty}\f1{n(n+1)}=0,$$
we deduce that
$$\lim_{n\to\infty}\log\f{\root{n+1}\of{u_{n+1}}}{\root n\of{u_n}}=0,\ \ \mbox{i.e.,}
\ \lim_{n\to\infty}\f{\root{n+1}\of{u_{n+1}}}{\root n\of{u_n}}=1.$$

For any $n\in\Z^+$, clearly
\begin{align*}\f{\root{n+1}\of{u_{n+1}}}{\root n\of{u_n}}>\f{\root{n+2}\of{u_{n+2}}}{\root{n+1}\of{u_{n+1}}}
&\iff\f{\log u_{n+1}}{n+1}-\f{\log u_n}n>\f{\log u_{n+2}}{n+2}-\f{\log u_{n+1}}{n+1}
\\&\iff \Delta_n:=\f{2\log u_{n+1}}{n+1}-\f{\log u_n}n-\f{\log u_{n+2}}{n+2}>0.
\end{align*}
Observe that
\begin{align*}\Delta_n=&2\log\al+\f{2\log(1-\gamma^{n+1})}{n+1}-\f{2\log{\sqrt{\Delta}}}{n+1}
\\&-\l(2\log\al+\f{\log(1-\gamma^n)}n+\f{\log(1-\gamma^{n+2})}{n+2}-\f{\log\sqrt{\Delta}}n-\f{\log\sqrt{\Delta}}{n+2}\r)
\\=&\f{\log\Delta}{n(n+1)(n+2)}+\f2{n+1}\log(1-\gamma^{n+1})-\f{\log(1-\gamma^n)}n-\f{\log(1-\gamma^{n+2})}{n+2}.
\end{align*}
The function $f(x)=\log (1+x)$  on the interval $(-1,+\infty)$ is concave since $f''(x)=-1/(x+1)^2<0$. Note that $|\gamma|<1$.
If $-|\gamma|\ls x\ls0$, then $t=-x/|\gamma|\in[0,1]$ and hence
$$f(x)=f(t(-|\gamma|)+(1-t)0)\gs tf(-|\gamma|)+(1-t)f(0)=qx, $$
where $q=-\log(1-|\gamma|)/|\gamma|>0$. Note also that $\log(1+x)<x$ for $x>0$. So we have

\begin{align*}
& \log \left(1- \gamma^{n+1} \right) \gs \log \left(1- |\gamma|^{n+1} \right) \gs - q|\gamma|^{n+1}, \\
& \log \left( 1- \gamma^{n} \right) \ls \log \left( 1 + |\gamma|^{n} \right) < |\gamma|^{n}, \\
& \log \left(1- \gamma^{n+2} \right) \ls \log \left( 1 + |\gamma|^{n+2} \right) < |\gamma|^{n+2}.
\end{align*}
Therefore
\[
\Delta_n > \frac{\log \Delta}{n(n+1)(n+2)} -|\gamma|^{n} \left( \f{2q|\gamma|}{n+1} + \f1n + \f{|\gamma|^2}{n+2} \right)
\]
and hence
\begin{multline}\label{Dn-1.1}
n(n+1)(n+2)\Delta_n \\
> \log\Delta-|\gamma|^n\l(2q|\gamma|n(n+2)+(n+1)(n+2)+|\gamma|^2n(n+1)\r).
\end{multline}
Since $\lim_{n \to \infty} n^2 |\gamma|^n = 0$, when $\Delta>1$ we have $\Delta_n > 0$ for large $n$.

Now it remains to consider the case $\Delta=1$. Clearly $\gamma=(A-1)/(A+1)>0$.
Recall that $\log(1-x)<-x$ for $x\in(0,1)$.  As
$$\f{d}{dx}(\log(1-x)+x+x^2)=-\f1{1-x}+1+2x=\f{x(1-2x)}{1-x}>0\quad\t{for}\ x\in(0,0.5),$$
we have $\log(1-x)+x+x^2>\log1+0+0^2=0$ for $x\in(0,0.5)$.
If $n$ is large enough, then $\gamma^n<0.5$ and hence
$$\Delta_n=\f2{n+1}\log(1-\gamma^{n+1})-\f{\log(1-\gamma^n)}n-\f{\log(1-\gamma^{n+2})}{n+2}>w_n,$$
where
$$w_n:=\f 2{n+1}(-\gamma^{n+1}-\gamma^{2n+2})+\f{\gamma^n}n+\f{\gamma^{n+2}}{n+2}.$$
Note that
$$\lim_{n\to\infty}\f{nw_n}{\gamma^n}=-2\gamma+1+\gamma^2=(1-\gamma)^2>0.$$
So, for sufficiently large $n$ we have $\Delta_n>w_n>0$.

Now we show that $n\gs4$ suffices in the case $A=1$ and $B=-1$. Note
that $\Delta=5$ and $\gamma \approx -0.382$. The sequence
$(|\gamma|^n (n+1)(n+2))_{n \gs 1}$ is decreasing since
\[
|\gamma| \frac{(n+2)(n+3)}{(n+1)(n+2)} < \frac{1}{2} \left(1+\frac{2}{n+1} \right) \ls 1
\]
for $n \gs 1$.  It follows that $|\gamma|^n (n+1)(n+2) \ls
\gamma^6\times 7\times 8<1/3$ for $n \gs 6$. In view of
\eqref{Dn-1.1}, if $n\gs 6$ then
\begin{align*}
n(n+1)(n+2)\Delta_n >& \log 5 - |\gamma|^n (n+1)(n+2) \l(2q|\gamma|+1+|\gamma|^2\r)
\\>&\log5-\f{1+1+\gamma^2}3>\log 5 - 1 >0.
\end{align*}
It is easy to verify that $\Delta_4$ and $\Delta_5$ are positive. So
 $(\sqrt[n+1]{F_{n+1}} / \sqrt[n]{F_{n}})_{n \gs 4}$ is
strictly decreasing.

In view of the above, we have completed the proof of Theorem 1.1. \qed

\section{Proof of Theorem 1.2}

\noindent{\it Proof of Theorem 1.2}. Let $n\gs3$.
    It is well known that $|D_n-n!/e|\ls1/2$ (cf. \cite[p.\,67]{St}).
  Applying the Intermediate Value Theorem in calculus, we obtain
    \begin{align*}
        \left| \log D_n - \log\left( \frac{n!}{e} \right) \right| \ls \left|D_n-\frac{n!}{e}\right| \ls 0.5.
    \end{align*}
   Set $R_0(n)=\log D_n-\log n!$. Then
        $|R_0(n)| \ls 1.5.$

    Since $\lim_{n\to\infty}R_0(n)/n=0$, we have
    \begin{align*}\lim_{n\to\infty}\l(\f{\log D_{n+1}}{n+1}-\f{\log D_n}n\r)
    =&\lim_{n\to\infty}\l(\f{\log (n+1)!}{n+1}-\f{\log n!}n\r)
    \\=&\lim_{n\to\infty}\f{n\log(n+1)+n\log n!-(n+1)\log n!}{n(n+1)}
    \\=&\lim_{n\to\infty}\f{n\log n+n\log(1+1/n)-\log n!}{n(n+1)}
    \\=&\lim_{n\to\infty}\f{\log (n^n/n!)}{n(n+1)}.
    \end{align*}
    As $n!\sim\sqrt{2\pi n}(n/e)^n$ (i.e., $\lim_{n\to\infty}n!/(\sqrt{2\pi n}(n/e)^n)=1$) by Stirling's formula, we have
    $\log(n^n/n!)\sim n$ and hence
     $$\lim_{n\to\infty}\l(\f{\log D_{n+1}}{n+1}-\f{\log D_n}n\r)=0.$$
   Thus $\lim_{n\to\infty}\root{n+1}\of {D_{n+1}}/\root n\of{D_n}=1$.

   From the known identity $D_n/n!=\sum_{k=0}^n(-1)^k/k!$, we have the recurrence $D_n=nD_{n-1}+(-1)^{n}$ for $n>1$.
 Thus, if $n\gs3$ then
 $$R_0(n)-R_0(n-1)=\log\f{D_n}{n!}-\log\f{D_{n-1}}{(n-1)!}=\log\f{D_n}{nD_{n-1}}=\log\l(1+\f{(-1)^n}{nD_{n-1}}\r).$$

Fix $n\gs4$. If $n$ is even, then
$$0<R_0(n)-R_0(n-1)=\log\l(1+\f1{nD_{n-1}}\r)<\f1{nD_{n-1}}=\f1{D_n-1}\ls\f3{D_n+0.5}.$$
If $n$ is odd, then
$$0>R_0(n)-R_0(n-1)=\log\l(1-\f1{nD_{n-1}}\r)>\f{-2}{nD_{n-1}}=\f{-2}{D_n+1}\gs\f{-3}{D_n+0.5}$$
since $\log(1-x)+2x>0$ for $x\in(0,0.5)$. So
$$|R_0(n)-R_0(n-1)|<\f3{D_n+0.5}\ls \f{3e}{n!}$$
and hence
$$\bg|\f{R_0(n-1)-R_0(n)}{n-1}\bg|<\f{3e}{n!(n-1)}\ls\f{3e}{n(n-1)(n+1)}.$$
Similarly, we also have
$$\bg|\f{R_0(n+1)-R_0(n)}{n+1}\bg|<\f{3e}{n!(n+1)}\ls\f{3e}{n(n-1)(n+1)}.$$
Therefore,
\begin{align*}&\bg|\f{R_0(n+1)}{n+1}-\f{2R_0(n)}n+\f{R_0(n-1)}{n-1}-\f{2R_0(n)}{n(n-1)(n+1)}\bg|
\\=&\bg|\f{R_0(n+1)-R_0(n)}{n+1}+\f{R_0(n-1)-R_0(n)}{n-1}\bg|\ls\f{6e}{n(n-1)(n+1)}
\end{align*}
and hence
$$\bg|\f{R_0(n+1)}{n+1}-\f{2R_0(n)}n+\f{R_0(n-1)}{n-1}\bg|\ls\f{2|R_0(n)|+6e}{n(n-1)(n+1)}\ls\f{6e+3}{n(n-1)(n+1)}.$$
Thus $|R_1(n)|\ls 6e+3$, where
$$R_1(n):=n(n-1)(n+1)\l(\f{R_0(n+1)}{n+1}-\f{2R_0(n)}n+\f{R_0(n-1)}{n-1}\r).$$

    Since
    \begin{align*}
        \log( (n-1)! )=&\sum_{k=1}^{n-1}\int_k^{k+1}(\log k) dx
        \\<&\sum_{k=1}^{n-1}\int_k^{k+1}\log xdx
        =\int_1^n \log x dx = n \log n - n + 1
        \\<&\sum_{k=1}^{n-1}\int_k^{k+1}(\log(k+1))dx=\log(n!),
    \end{align*}
    we have
    \begin{align*}
        n \log n - n < \log(n!) = \log( (n-1)! ) + \log n < n \log n - n + \log n + 1
    \end{align*}
    and so
  $\log(n!) = n \log n - n + R_2(n)$
    with $|R_2(n)| < \log n + 1$.

    Observe that
    \begin{align*}
        &\quad \frac{\log D_{n+1}}{n+1}  - \frac{2}{n} \log D_n + \frac{\log D_{n-1}}{n-1} \\
        &= \frac{\log (n+1)!}{n+1}  - \frac{2\log n! }{n} + \frac{ \log (n-1)!}{n-1} + \frac{R_1(n)}{(n-1) n (n+1)}\\
    &= \frac{2\log n!}{(n-1) n (n+1)}  - \frac{\log n }{n-1} + \frac{\log(n+1)}{n+1} + \frac{R_1(n)}{(n-1) n (n+1)}\\
    &= - \frac{2n}{(n-1)n(n+1)} + \frac{\log (n+1) - \log(n)}{n+1}
    +\frac{2 R_2(n)+R_1(n)}{(n-1) n (n+1)} \\
    &\ls - \frac{2n}{(n-1)n(n+1)} + \frac{n-1}{(n-1)n(n+1)}
    + \frac{2 R_2(n)+R_1(n)}{(n-1) n (n+1)}
    \\ &= - \frac{n + 1 - 2 R_2(n)- R_1(n)}{(n-1)n(n+1)}.
    \end{align*}
    If $n \gs 27$, then $n + 1 - 2 R_2(n)- R_1(n) > n - 2 \log n - 1- 6e -3 > 0$, and hence we get
      $$\log\f{\root{n}\of{D_{n}}}{\root{n-1}\of{D_{n-1}}}>\log\f{\root{n+1}\of{D_{n+1}}}{\root{n}\of{D_{n}}}.$$
    By a direct check via computer, the last inequality also holds for $n=4,\ldots,26$. Therefore,
 the sequence $(\sqrt[n+1]{D_{n+1}}/\sqrt[n]{D_n})_{n\gs3}$ is strictly decreasing.
This ends the proof. \qed

\section{Proof of Theorem 1.3}

\begin{lemma} \label{Lem4.1} For $x>0$ we have
\begin{equation}\log (1+x)>x-\f{x^2}2.\label{4.1}\end{equation}
\end{lemma}
\Proof. As
$$\f d{dx}\l(\log (1+x)-x+\f{x^2}2\r)=\f{x^2}{1+x},$$
we see that $\log(1+x)-x+x^2/2>\log1-0+0^2/2=0$ for any $x>0$. \qed

\begin{lemma} \label{Lem4.2} Let $m,n\in\Z^+$ with $n\gs3$. If $m\gs11$ or $n\gs30$, then
\begin{equation}H_n^{(m)}\log H_n^{(m)}>4\l(\f 2{n+2}\r)^{m-1}.\label{4.2}\end{equation}
\end{lemma}
\Proof. Recall that $H_n$ refers to $H_n^{(1)}$. If $n\gs30$, then
$H_n\log H_n\gs H_{30}\log H_{30}>4$
and hence (\ref{4.2}) holds for $m=1$.

Below we assume that $m\gs2$.
As $n\gs3$, we have $H_n^{(m)}\log H_n^{(m)}\gs H_3^{(m)}\log
H_3^{(m)}.$ So it suffices to show that
\begin{equation}\l(\f{n+2}2\r)^{m-1}H_3^{(m)}\log H_3^{(m)}>4\label{4.3}\end{equation}
whenever $m\gs11$ or $n\gs30$.
By Lemma 4.1,
\begin{align*}\log H_3^{(m)}=\log(1+2^{-m}+3^{-m})>&2^{-m}+3^{-m}-\f{(2^{-m}+3^{-m})^2}2
\\>&2^{-m}+3^{-m}-\f{(2^{1-m})^2}2=\f1{2^m}+\f1{3^m}-\f 2{4^m}.
\end{align*}
If $m\gs3$, then $(4/3)^m\gs(4/3)^3>2$ and hence $\log
H_3^{(m)}>1/2^m$. Note also that $H_3^{(2)}\log H_3^{(2)}>1/4$. So
we always have
$$H_3^{(m)}\log H_3^{(m)}>\f1{2^m}.$$
If $m\gs11$, then $1.25^m\gs1.25^{11}>10$ and hence
$$\f 1{2^m}>\f{4}{2.5^{m-1}}\gs\f 4{((n+2)/2)^{m-1}},$$
therefore (\ref{4.3}) holds. When $n\gs 30$, we have
$$\f1{2^m}\gs\f{1}{2^{4m-6}}=\f4{16^{m-1}}\gs 4\l(\f2{n+2}\r)^{m-1}$$
and hence (\ref{4.3}) also holds. \qed

\medskip
\noindent{\it Proof of Theorem 1.3}. Let $m\gs1$ and $n\gs3$. Set
$$\Delta_n(m) :=\log\f{\root{n+1}\of{H_{n+1}^{(m)}}}{\root{n}\of{H_{n}^{(m)}}}-\log\f{\root{n+2}\of{H_{n+2}^{(m)}}}{\root{n+1}\of{H_{n+1}^{(m)}}}
=\frac{2}{n+1} \log H_{n+1}^{(m)} - \frac{ \log H_n^{(m)}}{n} - \frac{\log H_{n+2}^{(m)}}{n+2}.$$
It suffices to show that $\Delta_n(m)<0$. This can be easily verified by computer if $m\in\{1,\ldots,10\}$ and $n\in\{3,\ldots,29\}$.

Below we assume that $m\gs11$ or $n\gs 30$.
Recall (\ref{4.1}) and the known fact that  $\log (1+x)<x$ for $x>0$.
We clearly have
$$\log\f{H_{n+1}^{(m)}}{H_n^{(m)}}=\log \left( 1 + \frac{1}{(n+1)^m H_n^{(m)}} \right) < \frac{1}{(n+1)^m H_n^{(m)}}$$
and
$$\log\f{H_{n+2}^{(m)}}{H_n^{(m)}}>\log \left(1 + \frac{2}{(n+2)^m H_n^{(m)}} \right) > \frac{2}{(n+2)^m H_n^{(m)}} - \frac{2}{(n+2)^{2m} (H^{(m)}_n)^2}.$$
It follows that
\begin{align*}
\Delta_n(m)  =& \left( \frac{2}{n+1} -  \frac{1}{n} -  \frac{1}{n+2} \right) \log H_n^{(m)} \\
& + \frac{2}{n+1} \log \f{H_{n+1}^{(m)}}{H_n^{(m)}} - \frac{1}{n+2} \log \f{H_{n+2}^{(m)}}{H_n^{(m)}} \\
 <& \frac{-2 \log H_n^{(m)}}{n(n+1)(n+2)} + \frac{2}{(n+1)^{m+1} H_n^{(m)}}
\\& - \frac{2}{(n+2)^{m+1} H_n^{(m)}} + \frac{2}{(n+2)^{2m+1} (H_n^{(m)})^2}.
\end{align*}
Since $(n+2)^{m+1}=\sum_{k=0}^{m+1}\bi{m+1}k(n+1)^k$ by the binomial theorem, we obtain
\begin{align*}
\Delta_n(m) \ls& \frac{-2 \log H_n^{(m)}}{n(n+1)(n+2)} + \frac{2\sum_{k=0}^m\bi{m+1}k(n+1)^k}{(n+1)^{m+1} (n+2)^{m+1} H_n^{(m)}}
+ \frac{2}{(n+2)^{m+2}H_n^{(m)}} \\
 <&  \frac{-2 \log H_n^{(m)}}{n(n+1)(n+2)} + \frac{2(n+1)^m\sum_{k=0}^m\bi{m+1}k}{(n+1)^{m+1}(n+2)^{m+1} H_n^{(m)}}
 + \frac{2}{(n+1)(n+2)^{m+1} H_n^{(m)}}
 \\=&\frac{-2 \log H_n^{(m)}}{n(n+1)(n+2)} + \frac{2(2^{m+1}-1)+2}{(n+1)(n+2)^{m+1} H_n^{(m)}}.
\end{align*}
Thus
\begin{align*}n(n+1)(n+2)\Delta_n(m)\f{H_n^{(m)}}2<&-H_n^{(m)}\log H_n^{(m)}+\f{2^{m+1}n}{(n+2)^m}
\\<&4\l(\f2{n+2}\r)^{m-1}-H_n^{(m)}\log H_n^{(m)}.
\end{align*}
Applying (3) we find that $\Delta_n(m)<0$.
This completes the proof. \qed
\medskip

\noindent{\bf Acknowledgments}. The initial version of this paper was posted to arXiv in 2012 as a preprint with the ID arXiv:1208.3903.
The authors are grateful to the two referees for their helpful comments.

\end{document}